\documentclass{article}
\usepackage{mathmag}

\usepackage{amsmath,amsthm}
\usepackage{graphicx}
\usepackage{hyperref}
\usepackage{url}
\usepackage{amsfonts}

\usepackage{multicol}

\theoremstyle{theorem}

\theoremstyle{definition}

\allowdisplaybreaks

\makeatletter
\@addtoreset{footnote}{page}
\makeatother

\begin{document}

\title{Cooking Poisons: Thinking Laterally with Game Theory}

\author{Timothy Y. Chow\\
\scriptsize Center for Communications Research\\
Princeton, NJ 08540\\
tchow@alum.mit.edu}

\maketitle


\section{Rabin's lateral-thinking puzzle}

In the late 1980s, a lateral-thinking puzzle by Michael Rabin
was posted to a Carnegie Mellon University electronic bulletin board.
The puzzle
seems to have been largely forgotten, but deserves to be better known.
Paraphrased, it went something like this.

There is a world in which the inhabitants have a strange physiology. A
healthy person who ingests a poison will die within an hour unless he
or she ingests a stronger poison; the stronger poison restores complete
health.
There are two distinct types of poisons: magical and medical,
which are dispensed by the Royal Magician and the Royal Physician
respectively (no poison is both magical and medical).
Moreover, all poisons in this world are strictly linearly ordered in strength.
All of these facts are common knowledge.

The King decides that he wants to find the strongest poison in the land,
because it will not only be very useful for eliminating enemies but will
also act as an antidote to any other poison. So he calls in the Royal
Magician and the Royal Physician and says, ``I want each of you to return
here to my royal chambers at noon one week from now. Bring a vial of your
strongest poison. To give you incentive to bring your strongest poison,
I will require each of you to drink from the other's vial
first, and then drink from your own vial. I will have trained observers
present to make sure that you cannot cheat. Then you will be watched for
one hour, during which you may not ingest any substances. The person
who has the stronger poison will of course survive, and the other will
die. The death of the loser is unfortunate, but I have decided that
it is a necessary sacrifice to ensure national security.
If I detect any attempt to circumvent these rules,
then you will both be executed.
You may go now, but you must return at the specified time.''

The Royal Physician and the Royal Magician leave, both very disturbed,
because neither wants to die. Each has had some experience with the
other's poisons and knows that some of them are quite potent. Neither
one is fully confident of having the strongest poison. Nor does
either have any way of getting access to the other's poisons. They rack
their brains all week trying to think how they can best ensure their
own survival.

The appointed time arrives, and the two Royal Servants return. They follow
the specified protocol exactly, and are watched carefully for one hour. To
everyone's astonishment, both Royal Servants keel over and die within
the hour. The Royal Coroner confirms that both died of poisoning.

\textit{What happened?}

\section{What is a lateral-thinking puzzle?}

Before reading on, you are encouraged to try to solve the puzzle yourself!

A \emph{lateral-thinking puzzle} is a puzzle whose solution cannot be
arrived at by straightforward logical or mathematical reasoning, but
requires creative thinking.
The term \emph{lateral thinking} was popularized by
psychologist Edward de Bono~\cite{debono} around 1971;
the interested reader can easily find many books full of
lateral-thinking puzzles, including some with a
mathematical focus~\cite{machale}.
A good lateral-thinking puzzle must strike
a delicate balance between spelling out too many details (which may
give away the solution) and leaving too many details unspecified
(which may open the door to unintended solutions that appear to solve
the puzzle as stated, but in an unsatisfying way).

For example, when I first encountered Rabin's puzzle, there was nothing
in the problem statement to rule out the possibility that the two Royal
Servants brought identical poisons.  This was not the intended solution,
so in my version, I have blocked it by
stipulating that no poison has exactly the same strength as another poison.

\section{The intended solution}

Each Royal Servant drank a weak poison just before arriving at the
showdown, and, instead of bringing a strong poison as requested, brought
water. They drank the other Servant's water, then drank their own water,
and died of their own poison.

But why did they do this?  The answer is that the Magician was hoping that
the Physician would not think of the same trick and would, as requested,
bring a strong poison to the showdown; that way, the Physician's strong
poison would cure the Magician of the Magician's own weak poison,
and the Magician would live,
whereas the Physician would die from drinking the Physician's own poison.
(There would remain
the risk that the King would discover the attempt to circumvent the rules,
but this risk could be mitigated
by drinking every last drop of water.)
The Physician thought similarly, and so both died.

I find this solution to be very satisfying, and I congratulate Michael
Rabin on crafting a fine puzzle.  But as we shall now see, the puzzle
turns out to be much richer than Rabin anticipated!

\section{Cooking the puzzle}

Chess problemists use the term \emph{cook} for an unintended solution
to a chess problem.  In the chess problem world, a cook is regarded as a fatal flaw
that invalidates the problem.  In the world of lateral thinking, however, a
cook may be a feature rather than a bug, if it is just as satisfying as the
intended solution.  We shall soon see that there are
indeed several alternative solutions to Rabin's puzzle that are in the same spirit
as the intended solution.

But first, let us be explicit about a point that was left somewhat
vague in the original narrative.  What happens if a healthy person drinks a
poison, and then drinks a second poison that is \emph{weaker} than the poison
that made them ill in the first place?  In what follows, we will make the
assumption that there are only two states that a (living) person can be in:
healthy, or mortally ill.  We further assume that if
a healthy person drinks Poison~X and becomes
mortally ill, then drinking poisons that are weaker than Poison~X,
or drinking further doses of Poison~X, leaves the person
in a mortally ill state, which will be cured if and only if
the person subsequently drinks a poison stronger than Poison~X.
For example, suppose Poison~X is stronger than Poison~Y, which in turn
is stronger than Poison~Z. Then a healthy person who drinks
Poison~Y, then Poison~Z, then Poison~X, ends up healthy,
but a healthy person who drinks Poison~Y, then Poison~X,
then Poison~Z, ends up mortally ill.
This assumption simplifies the analysis without (in my
opinion) ruining any essential features of the puzzle.

The key to unlocking the deeper secrets of Rabin's puzzle is what we might
humorously call \emph{Sicilian reasoning}, after a well-known scene in the 
movie \emph{The Princess Bride} in which Vizzini, a Sicilian, tries to
outsmart the hero by convoluted arguments of the form, ``I know
that you know that I would think that you would think that$\ldots\,$.''
Let us call the strategy of arriving healthy and bringing a poison
the \emph{Conventional} strategy (or Strategy~C for short).
In the intended solution of the puzzle, each Servant realizes that
there is an alternative strategy which defeats Strategy~C, namely
the \emph{Advanced} strategy (or Strategy~A for short) of
drinking a poison in advance and bringing water.
Alas, they both employ Strategy~A and they both die.

But we can take a tip from Vizzini, and go a step further.
I might reason as follows.  ``If my opponent employs Strategy~C,
then employing Strategy~A is a surer way to survive than
employing Strategy~C myself and gambling that I have the strongest poison.
But my opponent is no fool, and will probably realize the same thing,
and employ Strategy~A rather than Strategy~C.
I need to outsmart my opponent and find a way to defeat Strategy~A.''
Is there a strategy that will defeat Strategy~A?
Yes!  The \emph{Blank} strategy (or Strategy~B for short),
of bringing water but drinking \emph{nothing} in advance
(i.e., arriving at the showdown in full health), will defeat Strategy~A.
Unfortunately, Strategy~B is suicidal against Strategy~C,
but if I feel it is more likely than not that
my opponent will employ Strategy~A,
then I will prefer Strategy~B to Strategy A or~C,
which both lose against Strategy~A.

Vizzini would not stop here.  What if I am convinced
that my opponent will employ Strategy~B?  In that case,
if I myself were to employ Strategy~B,
then we would both survive the drinking ordeal, but then the King would
know that at least one of us had cheated, and would execute us both.
My only chance would be to employ the \emph{Double Dose} strategy
(or Strategy~D for short) of drinking a weak poison in advance,
and bringing a stronger poison to the showdown.
Strategy~D beats Strategy~B,
because my strong poison would cure me and kill my opponent.
Of course, if my opponent uses Strategy~A, then I should use
Strategy~B rather than Strategy~D$\ldots\,$.

Are there further strategies to consider beyond A, B, C, and~D?
No, not without introducing further creative possibilities beyond
those needed for the intended solution to the puzzle.
Each Servant may either drink a poison or drink nothing before the
showdown, and may either bring a poison or bring water to the showdown.
These four possibilities are precisely our four strategies
A, B, C, and~D.  Let us also point out that in Strategy~D,
bringing a poison equal to or weaker than the poison drunk in advance
can never be a good idea no matter what the opponent does,
so we always assume that in Strategy~D, a strictly weaker poison
is ingested ahead of time.

\section{The alternative solutions}

All four strategies have some merit, and could plausibly be chosen
by either Servant.
Going through all the possibilities one by one, we find that
there are three ways (other than the
intended ``A-versus-A'' scenario) in which both Servants
could keel over dead from poisoning.

\begin{enumerate}
\item One Servant employs Strategy~C and the other Servant employs Strategy~B.
\item One Servant employs Strategy~C and the other Servant employs Strategy~D,
and the strength of the former Servant's poison is in between the strengths
of the latter Servant's two poisons.
\item Both Servants employ Strategy~D, and each Servant's weaker poison is
weaker than the opponent's stronger poison.
\end{enumerate}

There is technically a fourth possibility, which is that one Servant
employs Strategy~A, the other Servant employs Strategy~C, and the
poison drunk in advance by the Servant employing Strategy~A is stronger than
the poison brought by the Servant employing Strategy~C.
However, we can argue that this possibility violates the spirit of the
puzzle, by putting together three observations.

\begin{enumerate}
\renewcommand{\theenumi}{\alph{enumi}}
\item \label{obs:a}
\emph{A Servant employing Strategy~A will use the weakest
available poison.} The reason is that the hope is to be cured
by the opponent's poison, so there is never an incentive to use
anything but the weakest poison.
\item \label{obs:b}
\emph{A Servant employing Strategy~C will use the strongest
available poison.} Doing so maximizes the chances
against Strategies C and~D,
and against Strategies A and~B,
Strategy~C is suicidal regardless.
\item \label{obs:c}
\emph{It is not the case that one Servant's poisons
are all stronger than all of the other Servant's poisons.}
Indeed, we are told
that each Servant has had some experience with the
other's poisons and knows that some of them are quite potent.
\end{enumerate}

\section{Game theory}

As far as finding cooks is concerned, our analysis is complete,
but there is a natural question that remains unanswered.
Do we have any advice to give the Servants,
or are they stuck with endlessly circular Sicilian reasoning,
with no reason to prefer one strategy over another?

Textbooks on game theory
(e.g., \cite{karlin-peres,maschler-solan-zamir})
tell us that if every \emph{deterministic} (or \emph{pure})
strategy available to me can be defeated by my opponent, then it may be 
advantageous for me to employ a \emph{randomized} (or \emph{mixed}) strategy.
We all know that in the game of Rock, Paper, Scissors (or Roshambo),
it makes no sense to always make the same play; it is essential to
randomize one's choices in order to keep the opponent guessing.

To apply the tools of game theory to Rabin's puzzle, we need to specify the
(expected) \emph{payoffs} to each Servant in each possible scenario
(A-versus-A, A-versus-B, etc.).  It seems reasonable
to assume that each Servant wishes to maximize their own chances of survival,
so let us assign a payoff of~1 to a Servant who survives and whose opponent dies,
and a payoff of~0 otherwise (so in particular, if both Servants die of poisoning,
or if both Servants survive the drinking ordeal and are therefore executed by
the King, then both Servants receive a payoff of~0).
Probability of survival is then equivalent to expected payoff,
which is what game theory is designed to maximize.

At this point, we run into a complication.
Under our assumptions, there are some scenarios in which
we can say for certain who will live and who will die,
regardless of the strengths of the poisons used (if any);
in particular, this is true whenever at least one Servant employs
either Strategy~A or Strategy~B. But in other scenarios,
such as C-versus-C or D-versus-D, the outcome depends on the relative
strengths of the Servants' poisons.  How do we deal with
the uncertainty surrounding the relative strengths of the various
available poisons?

A full answer is very complicated, so let us analyze a special case.
Let us assume that both Servants have exactly two poisons
at their disposal.  If we denote the Magician's poisons with an~M
and the Physician's poisons with a~P, then there are four possible
strength orderings of these four poisons:
MPMP (i.e., the Magician's stronger poison is the strongest,
followed by the Physician's stronger poison, followed by the
Magician's weaker poison, followed by the Physician's weaker
poison in last place), MPPM, PMPM, and PMMP
(per Observation~\ref{obs:c} in the previous section,
we rule out MMPP and PPMM).
Given one of these four rank orderings, we can write down the
payoff matrices for each Servant.
It turns out that the payoff matrices for MPPM
are the same as for MPMP;
see Figure~\ref{fig:MPPM},
where the rows give the Magician's strategy
and the columns give the Physician's strategy.
For example, the 1 in row~A and column~C
of the matrix on the left means that if the Magician chooses Strategy~A
and the Physician chooses Strategy~C, then the Magician will survive
and the Physician will die.

\begin{figure}[ht]
\begin{minipage}{.5\linewidth}
   \centering
   \begin{tabular}{c | c c c c}
   \ & A & B & C & D\\
   \hline
   A & 0 & 0 & 1 & 0$\vphantom{X^{X^X}}$\\
   B & 1 & 0 & 0 & 0\\
   C & 0 & 0 & 1 & 1\\
   D & 0 & 1 & 0 & 0
   \end{tabular}
\end{minipage}%
\begin{minipage}{.5\linewidth}
   \centering
   \begin{tabular}{c | c c c c}
   \ & A & B & C & D\\
   \hline
   A & 0 & 1 & 0 & 0$\vphantom{X^{X^X}}$\\
   B & 0 & 0 & 0 & 1\\
   C & 1 & 0 & 0 & 0\\
   D & 0 & 0 & 0 & 0
   \end{tabular}
\end{minipage}
\caption{MPPM or MPMP payoff matrices (Magician on left, Physician on right)}
\label{fig:MPPM}
\end{figure}

To obtain the payoff matrices for PMMP or PMPM, we simply swap the two
matrices while transposing rows and columns.  See Figure~\ref{fig:PMMP},
where again the rows give the Magician's strategy and the columns give the
Physician's strategy.

\begin{figure}[ht]
\begin{minipage}{.5\linewidth}
   \centering
   \begin{tabular}{c | c c c c}
   \ & A & B & C & D\\
   \hline
   A & 0 & 0 & 1 & 0$\vphantom{X^{X^X}}$\\
   B & 1 & 0 & 0 & 0\\
   C & 0 & 0 & 0 & 0\\
   D & 0 & 1 & 0 & 0
   \end{tabular}
\end{minipage}%
\begin{minipage}{.5\linewidth}
   \centering
   \begin{tabular}{c | c c c c}
   \ & A & B & C & D\\
   \hline
   A & 0 & 1 & 0 & 0$\vphantom{X^{X^X}}$\\
   B & 0 & 0 & 0 & 1\\
   C & 1 & 0 & 1 & 0\\
   D & 0 & 0 & 1 & 0
   \end{tabular}
\end{minipage}
\caption{PMMP or PMPM payoff matrices (Magician on left, Physician on right)}
\label{fig:PMMP}
\end{figure}

We make one more assumption, that the Magician and Physician
are equally likely to have the strongest poison.
So the matrices of expected payoffs are obtained by averaging
the payoff matrices in Figures \ref{fig:MPPM} and~\ref{fig:PMMP}.
See Figure~\ref{fig:overall}.

\begin{figure}[ht]
\begin{minipage}{.5\linewidth}
   \centering
   \begin{tabular}{c | c c c c}
   \ & A & B & C & D\\
   \hline
   A & 0 & 0 & 1 & 0$\vphantom{X^{X^X}}$\\
   B & 1 & 0 & 0 & 0\\
   C & 0 & 0 & 1/2 & 1/2 \\
   D & 0 & 1 & 0 & 0
   \end{tabular}
\end{minipage}%
\begin{minipage}{.5\linewidth}
   \centering
   \begin{tabular}{c | c c c c}
   \ & A & B & C & D\\
   \hline
   A & 0 & 1 & 0 & 0$\vphantom{X^{X^X}}$\\
   B & 0 & 0 & 0 & 1\\
   C & 1 & 0 & 1/2 & 0\\
   D & 0 & 0 & 1/2 & 0
   \end{tabular}
\end{minipage}
\caption{Overall payoff matrices (Magician on left, Physician on right)}
\label{fig:overall}
\end{figure}

\section{Nash equilibria}

All right, we have the payoff matrices.  What is each Servant's optimal strategy?

In a \emph{zero sum game}, one player's payoff is the negative of the
other player's payoff (possibly after a constant is added to all payoffs).
Zero sum games admit a fairly clear notion of an optimal strategy,
but unfortunately, our game is not a zero sum game, because sometimes
both players receive a payoff of~0, while other times one player
receives a payoff of~1 and the other receives a payoff of~0.
The meaning of an ``optimal'' strategy---if there is such a thing---is
therefore not immediately clear.

Game theorists often suggest
that a necessary condition for a strategy---or really, a pair
of strategies, one for each player---to be optimal is that it
be a \emph{Nash equilibrium}.  In a Nash equilibrium, by definition,
neither player can increase their expected payoff by changing strategies
while the opponent's strategy remains fixed. The nice thing about the
Nash equilibrium concept is that, under very mild conditions (which are
certainly satisfied in our game), at least one Nash equilibrium always
exists.

Our game is small enough that the
Nash equilibria can be computed by hand from 
first principles, though doing so is rather laborious.
There is software for computing Nash equilibria;
e.g., Gambit (\texttt{http://www.gambit-project.org}),
or Game Theory Explorer (\texttt{http://app.test.logos.bg}).
However we compute them,
we find that there are exactly three Nash equilibria.

\begin{enumerate}
\item Choose one of the four strategies uniformly at random.
\item The Magician chooses either B or~C with probability~1/2,
and the Physician chooses A with probability~1/3 and
D with probability~2/3.
\item The Physician chooses either B or~C with probability~1/2,
and the Magician chooses A with probability~1/3 and
D with probability~2/3.
\end{enumerate}

It is easy to check that in the first equilibrium, each Servant has
a survival probability of~1/4.
The other two equilibria are more interesting.
The A/D Servant has a survival probability of~1/2
and the B/C Servant has a survival probability of~1/3.

So if both Servants get as far as calculating the Nash equilibria,
what should they do?  Unfortunately, game theory does not
always give us a \emph{prescription} for how to maximize our
expected payoff, especially if there are multiple Nash equilibria,
as is the case here.
At first glance, it might seem that both Servants would prefer
one of the latter two equilibria to the first, since both of
their survival probabilities would go up.
But the situation is an adversarial one,
and they cannot unilaterally force their way from
one equilibrium to another, even if both would benefit from the switch.

To make matters even more complicated,
an argument can even be made for the \emph{non-equilibrium} strategy
that chooses B, C, D with probabilities 1/4, 1/2, 1/4 respectively,
because it guarantees a survival probability of 1/4 no matter
what the opponent does (and there is no strategy that can \emph{guarantee}
a survival probability larger than~1/4).  But it is hard to argue
that it is ``optimal'' for both Servants to adopt this strategy,
since one Servant could then unilaterally increase survival
probability by always playing~A.

\section{Concluding remarks}

It may feel unsatisfying to end this article without
a definitive recommendation for the Servants to follow,
but I regard this as a good thing, because it creates
an opportunity for you, the reader, to do further research.
We have not exhausted all the ramifications of Rabin's puzzle,
and perhaps further investigation will shed more light on
what the Servants ought to do.

In particular, new subtleties arise when more poisons are available.
Recall that in our simplified example, 
we ruled out the possibility that one Servant's two poisons were
both stronger than the other Servant's two poisons,
and so D-versus-D would always kill both Servants.
But if the Servants have access to more than two poisons,
then it is not implausible that one Servant might have
the two strongest poisons in the kingdom,
leading to a victory in a D-versus-D showdown.
Furthermore, the decision about \emph{which} two poisons
to select when employing Strategy~D is not obvious.
Against Strategy~C, it is best to use one's two strongest poisons,
but against Strategy~A, it is best to use one's two
\emph{weakest} poisons.

Empirical investigation of the game would also be interesting.
What happens in practice when experienced players repeatedly
play against each other (perhaps re-randomizing the strengths
of the poisons from one round to the next)?
An incredibly fruitful stream of research in game theory
was initiated in 1980 when Robert Axelrod~\cite{axelrod} set up a
tournament for the iterated prisoner's dilemma.
Maybe an enterprising reader will be motivated to
set up a tournament, and invite participants
to \emph{choose their poison!}

\section{Acknowledgments}

I would like to thank the anonymous referee for numerous suggestions
which improved the exposition of this paper.

\begin{abstract}
We revive an old lateral-thinking puzzle by Michael Rabin,
involving poisons with strange properties.
We show that the puzzle admits several unintended solutions
that are just as interesting as the intended solution.
Analyzing these alternative solutions using game theory
yields surprisingly subtle results and several unanswered questions.
\end{abstract}

\begin{biog}
\item[Timothy Y. Chow] (MR Author ID 336506)
received his Ph.D. from MIT in 1995 under
Richard Stanley and currently works at the Center for Communications
Research in Princeton.  His research interests include algebraic
combinatorics and computational complexity theory.  He would
publish more if he didn't spend far too much time playing
backgammon.
\end{biog}
\vfill\eject

\end{document}